\documentclass[11pt]{article}
\usepackage{amsmath,amsthm,amssymb,latexsym,graphicx,amsfonts,amscd}

\usepackage{color}

\newtheorem{theorem}{Theorem}[section]

\newtheorem{problem}[theorem]{Problem}

\newtheorem{lemma}[theorem]{Lemma}

\numberwithin{equation}{section}

\def\Fp{\vadjust{}\penalty200 \hfill
\lower.3333ex\hbox{\vbox{\hrule\hbox{\vrule\phantom{\vrule height
6.83333pt depth 1.94444pt width 8.77777pt}\vrule}\hrule}}
\ifmmode\let\next\relax\else\let\next\par\fi \next}

\def\End{\mathop{\rm End}\nolimits}

\def\Hom{\mathop{\rm Hom}\nolimits}

\def\rk{\mathop{\rm rk}}

\def\N{{\mathbb N}}
\def\Z{{\mathbb Z}}
\def\Q{{\mathbb Q}}

\def\bS{{\mathbb S}}

\def\k{\kappa}

\def\l{\lambda}

\def\aln{{\aleph_0}}
\def\Cont{2^{\aln}}

\def\n+1d{{}^{ n+1 \downarrow }\l}

\def\size#1{\left|\,#1\,\right|}

\def\Rhat{\widehat R}

\def\to{\rightarrow}
\def\arr{\longrightarrow}

\def\Rhat{\widehat{R}}

\begin{document}

\title{\bf On cellular covers with free kernels}
\footnotetext{The first author was supported by the Spanish Ministry
of Education and Science MEC-FEDER grant MTM2007-63277.\\
The second author was supported by the project No. 963-98.6/2007
of the German-Israeli Foundation for Scientific Research \&
Development.

}

\author{Jos\'e L. Rodr\'{\i}guez and Lutz Str\"ungmann}
\date{January 7th, 2010}
\maketitle

\begin{abstract}
Recall that a homomorphism of $R$-modules $\pi: G\to H$ is called
a {\it cellular cover} over $H$ if $\pi$ induces an isomorphism
$\pi_*: \Hom_R(G,G)\cong \Hom_R(G,H),$ where $\pi_*(\varphi)= \pi
\varphi$ for each $\varphi \in \Hom_R(G,G)$ (where maps are
acting on the left). In this paper we show that every
cotorsion-free module $K$ of finite rank can be realized as the
kernel of a cellular cover of some cotorsion-free module of rank
$2$. In particular, every free abelian group of any finite rank
appears then as the kernel of a cellular cover of a
cotorsion-free abelian group of rank $2$. This situation is best possible in the sense that cotorsion-free abelian groups
of rank $1$ do not admit cellular covers with free kernel except
for the trivial ones. This work comes motivated by an example due
to Buckner and Dugas, and recent results obtained by
G\"obel--Rodr\'{\i}guez--Str\"ungmann, and Fuchs--G\"obel.
\end{abstract}

\section{Introduction}\label{introduction}
Recently cellular covers of groups and modules have attracted
much attention in the literature. Recall that a homomorphism
$\pi:G \to H$ of groups is a {\it cellular cover} over $H$, if
every homomorphism $\varphi: G\to H$ lifts uniquely to an
endomorphism $\widetilde\varphi$ of $G$ such that $\pi
\widetilde\varphi = \varphi$. In case that $\pi: G\to H$ is an
epimorphism we say that
$$
1\to K \to G\overset{\pi}{\to} H\to 1
$$
is a {\it cellular exact sequence}. Cellular covers are the
algebraic counterpart of cellular approximations of topological
spaces in the sense of J.H.C. Whitehead, or the more general
cellularization maps extensively studied in homotopy theory in
the 90s (see e.g. \cite{Bo97}, \cite{Cha96}, \cite{F97},
\cite{Hir02}, \cite{Nof99}). As for the dual case, namely for
localizations, there is sometimes a good interplay between
cellularization of spaces and cellularization of groups and
modules. This motivated a careful study of the algebraic setting
(see e.g. \cite{Cas00}, \cite{CRT98}, \cite{Flo},
\cite{RS01}, \cite{RS07}).\\
The general goal is to completely classify (up to isomorphism) all
possible cellular exact sequences for fixed cokernel $H$ or
kernel $K$. In some cases this is possible. For instance, if $H$
is divisible, then $G$ (above) can be determined explicitly as
was shown in \cite{CFGS07}. And if $H$ is torsion and reduced,
then the cellular exact sequence collapses and $K=0$, see
\cite{FG}.\\
However, if $H$ is reduced (and torsion-free), then $K$ must be
cotorsion-free, see \cite{BD}, \cite{FGSS07}, \cite{FG} and a
result by Buckner and Dugas \cite{BD} shows that any
cotorsion-free module $K$ is the kernel of arbitrarily large
cellular covers $G$ of size $\geq \Cont$. Dually, it was shown in
\cite{GRS09} that any cotorsion-free module $H$ of size $\geq
2^{\aleph_0}$ satisfying natural rigidity conditions is the
cokernel of arbitrarily large cellular covers. This makes a
classification
impossible.\\

Also the case of finite rank modules was treated in \cite{GRS09}
and also in \cite{FG}. If $R$ is a subgroup of the rational
numbers and $R$ is a ring, then $R$ does not admit cellular
covers except for the trivial ones (see \cite{FG}). However, if
$R$ is not a ring, then $R$ has arbitrarily large cellular covers
by \cite{FG}. Unfortunately, the proof in \cite{FG} of this result
contains a gap that we will fix in this paper (see Section $1$).
Moreover, the construction in \cite{FG} implies that the kernels of the cellular covers in this situation
are $p$-divisible for some prime $p$ and we will
prove that this is justified by showing the following (see
Theorem \ref{free-kernel}): If $R$ is a subgroup of $\Q$, then
$R$ does not admit cellular covers with free kernel.
Nevertheless, we prove in Section $3$ that any cotorsion-free
module $K$ of finite rank is in fact the kernel of a cellular
cover with cokernel of rank $2$ (see Theorem \ref{rank2}). That
the rank of $K$ has to be countable is a necessary restriction by
the following result. If $K$ is free, then $|K|\leq |H|$
(\cite{FGSS07}). This clarifies the situation for cellular covers
with cokernel of rank $1$ and $2$ almost completely since it is
left as an open question if also a free module of countably
infinite rank can be the kernel of some cellular cover of
cokernel of rank $2$.

\noindent

\section{Cellular covers of rational groups}\label{rank1}

In this section we consider a result from \cite{FG} and fill a gap
in its proof. Recall from \cite[Theorem 2.11]{FGSS07} that there
are torsion-free abelian groups $L$ that admit arbitrarily large
cellular covers \[ 0 \rightarrow K \rightarrow G \rightarrow L
\rightarrow 0.\] Necessarily, the kernel $K$ of the cellular cover
then has to be non-free if $|K|
> |L|$ by \cite[Proposition 1.4]{FGSS07}. This result was strengthened in \cite{FG} (to groups $L$ of rank one)
where the following was claimed. Recall that any rank one group
$S$ is a subgroup of the rational numbers $\Q$ and is of the form
$S=\left< \frac{1}{p^{n_p}} : n_p < m_p, p \in \Pi \right>$ for
some natural numbers or the $\infty$-symbol $m_p \in \N_0 \cup \{\infty\}$ where $\Pi$ is
the set of all primes and $\N_0$ is the set of natural numbers
including $0$. The type of $S$ is then the equivalence class (in
the sense of Baer) of the height sequence $(m_p-1 : p \in \Pi)$ (with the convention that $\infty -1=\infty$)
where two sequences $(m_p-1 : p \in \Pi)$ and $(k_p-1 :p \in
\Pi)$ are equivalent if the sum $\sum\limits_{p \in \Pi} m_p-k_p$
is finite. The types form a lattice that has a rich structure (see
\cite[Chapter XIII]{Fu}).

\begin{theorem}
\label{FGold} Let $R \subseteq \Q$ be a torsion-free abelian group
of rank 1.

\begin{itemize}
\item[(I)] Suppose that $R$ is not a ring. Then $R$ admits a cellular
cover of rank $\kappa$ for every cardinal $\kappa \geq 1$
(\cite[Lemma 5.3]{FG});
\item[(II)] Suppose $R$ is a ring. Then the only cellular covers of $R$
are the trivial ones (\cite[Theorem 6.1]{FG}). \end{itemize}
\end{theorem}

The idea of the proof for (I) is the following (see \cite[Proof
of Lemma 5.3]{FG}). For simplicity let $R=\left< \frac{1}{p} :
q\not= p \in \Pi \right>$ for some fixed prime $q$. The general case follows with easy modifications. For a cardinal
$\kappa \geq 1$ choose a torsion-free group $K$ of rank $\kappa$
such that $\End(K)\cong\Z[\frac{1}{q}]$ and $K$ does not contain
any pure subgroup of rank one whose type is at least the type of
$R$. The authors then choose elements $b_p \in K$ for $p \not= q$
such that $b_p$ is not divisible by $p$ in $K$ and put $G:=\left<
K, \frac{a+b_p}{p}: q\not=p \in \Pi \right>$ and claim that \[
(*) \quad 0 \rightarrow K \rightarrow G
\overset{\pi}{\rightarrow} R \rightarrow 0 \] is a cellular cover
with the obvious map $\pi$ sending $a$ onto $1 \in R$ and $K$ to
$0$. One easily checks that \begin{itemize}
\item $\End(R)\cong \Z$ \item $K$ is fully invariant in $G$ as it
is the maximal $q$-divisible subgroup of $G$
\item $\Hom(K,R)=0$ by the $q$-divisibility of $K$. \end{itemize}
Moreover, the authors also claim that \begin{itemize}
\item $\Hom(G,K)=0$ \end{itemize} which then would imply that the
above sequence $(*)$ is indeed a
cellular cover by \cite[Lemma 3.5]{FGSS07}.\\
However, the condition on the elements $b_p$ is not strong enough
to force $\Hom(G,K)=0$ as the following Lemma shows.

\begin{lemma}
If $K=\Z[\frac{1}{q}]$ and $b_p=p-1 \in K$, then the sequence
$(*)$ is not a cellular cover of $R$. However, $b_p \not\in pK$
for every $q \not=p \in \Pi$.
\end{lemma}

\begin{proof}
It is easy to see that in this case the sequence $(*)$ splits
since $\frac{a+b_p}{p}=\frac{a+p-1}{p}\equiv\frac{a-1}{p}$ mod $K$
and so $G=\left< K,\frac{a-1}{p} \right> \cong K \oplus R$. Hence
$(*)$ cannot be a cellular cover.
\end{proof}

This shows that one has to choose the elements $b_p$ more
carefully; we now show how.

\begin{proof}(of Theorem \ref{FGold} (I))
As above let $R=\left< \frac{1}{p} : q\not= p \in \Pi \right>$
for some fixed prime $q$. For a cardinal $\kappa \geq 1$ choose a
torsion-free group $K$ of rank $\kappa$ such that
$\End(K)\cong\Z[\frac{1}{q}]$ and $K$ is homogeneous of type
$\Z[\frac{1}{q}]$. Such a group exists for every cardinal
$\kappa$ as was pointed out in \cite{FG} (see also \cite{GT}). We
now divide the set of primes $\Pi$ into two disjoint infinite
subsets $\Pi_1$ and $\Pi_2$ with $q \in \Pi_1$. Fix a bijection
$\sigma: \Z[\frac{1}{q}] \rightarrow \Pi_2$ such that $z \not\in
\sigma(z)\Z[\frac{1}{q}]$ for all $0\not= z \in \Z[\frac{1}{q}]$.
By the type condition on $K$ we may now pick elements
$b_{\sigma(z)} \in K$ such that $zb_{\sigma(z)} \not\in
\sigma(z)K$ for all $0 \not= z \in \Z[\frac{1}{q}]$. Let
\[G=\left<K, \frac{a}{p}, \frac{a+b_r}{r} : p \in \Pi_1, r \in
\Pi_2 \right>. \] As above it is easy to see that $K$ is fully
invariant in $G$ and that $\Hom(K,R)=0$. We claim that now also
$\Hom(G,K)=0$. Therefore assume $\varphi \in \Hom(G,K)$. Then
$\varphi \restriction_K = z \cdot id_K$ for some $z \in
\Z[\frac{1}{q}]$ by the full invariance of $K$. Moreover, $\varphi(a)=k \in K$. We obtain
\[ \varphi(\frac{a}{p}) = \frac{k}{p} \in K \]
for all $p \in \Pi_1$. But $K$ does not contain non-zero elements
of type other than $\Z[\frac{1}{q}]$, hence $\varphi(a)=0$. Now
\[ \varphi(\frac{a+b_r}{r})=\frac{zb_r}{r} \in K \]
for all $r \in \Pi_2$. However, by the choice of $\sigma$ we have
\[ \frac{zb_{\sigma(z)}}{\sigma(z)} \not\in K \]
a contradiction if $z \not= 0$. Thus $z=0$ and hence $\Hom(G,K)=0$
and the sequence
\[ 0 \rightarrow K \rightarrow G \rightarrow R \rightarrow 0\]
is a cellular cover by \cite[Lemma 3.5]{FGSS07}.
\end{proof}

\section{Cellular covers of rational groups with free kernel}\label{rank1-free-kernel}

In the previous section it was shown that there is a cellular
cover

\[ 0 \rightarrow K \rightarrow G \rightarrow R \rightarrow 0 \]
of a rational group $R$ such that the kernel $K$ is isomorphic to
the rational group $\Z[\frac{1}{q}]$ for some fixed prime $q$. We
now ask if one can get the same result replacing
$\Z[\frac{1}{q}]$ by the integers $\Z$ or more generally by a
free group and still have a cokernel of rank one. The following
result shows that this is impossible.

\begin{theorem} \label{free-kernel}
Let $H \subseteq \Q$ be a torsion-free abelian group of rank 1 and
$H_0$ its nucleus. Suppose that $H$ is not a ring. Then $H$ does
not admit any cellular cover with kernel a free $H_0$-module.
\end{theorem}

\begin{proof} Let $H$ be as stated. Then $H \not\cong H_0$ and any cellular cover of $H$ is an $H_0$-module by \cite{FGSS07}. Without loss of generality we may assume that
$H_0=\Z$ and that $H=\langle \frac{1}{p}: p\in \Pi \backslash
\{2\} \rangle$. The general case is obtained by simple
modification of our arguments and therefore left to the reader.

Let $F= \langle e_i : i\in I\rangle$ be a free-abelian group, with
$I$ a non-empty index-set, and suppose that there exists a
cellular exact sequence
\begin{equation}\label{cellular sequence}
0\to F \longrightarrow G \stackrel{\pi}{\longrightarrow} H \to 0,
\end{equation}
that is, $\Hom(G,F)=0$ and every homomorphism $\varphi: G \to H$
lifts to a (unique) endomorphism $\psi: G \to G$ such that $\pi
\psi= \varphi$.

If $G$ fits in (\ref{cellular sequence}) then there exist elements
$a\in G$ and $z_p\in F$ for all $p\in \Pi$ such that:
 $$G = \langle F, \frac{a+z_p}{p} : p\in \Pi\rangle$$
where $\pi(a)=1 \in H$ and  $\pi(\frac{a+z_p}{p}) = \frac{1}{p}
\in H$. We express each $z_p$ as a linear combination of the base
elements of $F$, say $z_p= \sum\limits_{i \in I} z_p^i e_i$, where
almost all the coefficients are equal to $0$. We can assume
without loss of generality that every non-trivial coefficient
$z_p^i$ satisfies $(z_p^i,p)=1$, otherwise we could just erase it,
as the term $\frac{z_p^i e_i}{p}$ would belong to $F$.

We can now give explicit generators of $\Hom(G,H)$. For every
$r,s \in \Z$, and $i\in I$ define the homomorphism:
$$
\varphi=\varphi_{s,r,i}: G \to H \quad {\rm by }\quad
\varphi(a)=r, \quad \varphi(e_j)=0 \quad \textit{\rm for } j\neq
i, \quad \varphi(e_i)=s;$$ Note that $\varphi$ is well defined
since $(r+z_p^i s)$ is an integer and therefore $(r+z_p^i s)\in p
H$, for all $p\in \Pi\backslash\{2\}$ and $i\in I$. By linearity we then have
$$\varphi\left(\frac{a+z_p}{p}\right)= \frac{r+z_p^i s}{p}.
$$
Now, if $q\in \Pi$ is a fixed prime such that $r+z^i_qs \in q \Z$,
then one can define
$$
\varphi'=\varphi^q_{r,s,i}: G \to H \quad {\rm  by } \quad
\varphi'(a)=\frac{r}{q}, \quad  \varphi'(e_j)=0 \quad \textit{\rm
for } j\neq i, \quad \varphi'(e_i)=\frac{s}{q}, $$ and therefore
$$\varphi'\left(\frac{a+z_p}{p}\right)=\frac{r+sz_p^i}{pq}.$$

Note that the condition $r+z^i_qs \in q \Z$ is needed since
otherwise $\frac{r+sz_q^i}{q^2} \not\in H$. By assumption we know
that (\ref{cellular sequence}) is a cellular sequence, hence there
exist unique endomorphisms of $G$, $\psi=\psi_{r,s,i}$ and
$\psi'=\psi^q_{r,s,i}$ such that $\pi\psi=\varphi$ and
$\pi\psi'=\varphi'$. Torsion-freeness of $H$ and uniqueness of
lifting imply $q \psi'=\psi$.

We want to say more about the action of $\psi$. Obviously,
$$\varphi_{r,s,i}=\varphi_{r,0,i} +
\varphi_{0,s,i}=r\varphi_{1,0,i} + s \varphi_{0,1,i}.$$ If $s=0$,
then the unique lifting of $\varphi_{r,0,i}$ is multiplication by
$r$ on $G$. Similarly, uniqueness implies that
$s\psi_{0,1,i}=\psi_{0,s,i}$. It is now easy to check that there
exist elements $h_i$, $h_a\in F$ such that $\psi$ is given by:
\begin{equation}
\label{psi}
\begin{array}{rlll}
\psi :& G   & \longrightarrow  & G                              \\
      & a    & \mapsto         & ra + sh_a                      \\
      & e_j  & \mapsto         & re_j+sh_j \mbox{ \quad for $j\neq i$}\\
      & e_i  & \mapsto         & re_i+sh_i +sa \\
      & \frac{a+z_p}{p} & \mapsto        & \psi(\frac{a+z_p}{p})
\end{array}
\end{equation}
Note that $q\psi'=\psi$ yields in particular
\begin{equation}\label{psiei}
\psi(e_i)= re_i+sh_i +sa \in q G
\end{equation}
Since $\psi(a+z_p)\in p G$, we conclude for all $p\in \Pi$
\begin{equation}\label{psiazp}
\begin{array}{rcl}
\psi(a+z_p) & = & ra+sh_a+ \sum\limits_{j\neq i}z_p^j(re_j + sh_j) +z_p^i (re_i+sh_i+sa)\\
&&\\[-.2cm]
& = & r(a+z_p) +s(h_a + \sum\limits_{j} z_p^jh_j + z_p^ia) \in p
G.
\end{array}
\end{equation}
Using that $(a+z_p) \in p G$, we then obtain
$$
s(h_a + \sum_{j} z_p^jh_j + z_p^ia) \in p G.
$$
Subtracting $sz_p^i(a+z_p) \in p G$ from the previous expression
we get
\begin{equation}\label{star}
s(h_a + \sum_{j} z_p^jh_j - z_p^iz_p) \in p G \cap F = pF
\end{equation}
since $H$ is torsion-free and hence $F$ is pure in $G$.

Now recall that $\psi=q\psi'$ for the fixed prime $q$ satisfying
$(r+sz_q^i) \in q\Z$. Hence $\psi(a+z_p) \in qG$ for all primes $p
\in \Pi$. From (\ref{psiazp}) and adding and subtracting $sz_q^i
(a+z_p)$ we obtain
$$
\begin{array}{rcl}
\psi(a+z_p) & =  & \psi(a+z_p) + sz_q^i (a+z_p) - sz_q^i (a+z_p) \\
            & =  & (r+sz_q^i)(a+z_p) + s[h_a + \sum\limits_{j} z_p^jh_j + z_p^ia] -
            sz_q^i            (a+z_p)\\
            & =  & (r+sz_q^i)(a+z_p) + s[h_a + \sum\limits_{j} z_p^jh_j + z_p^ia - z_q^ia - z_q^i z_p] \in qG. \\
\end{array}
$$
Since $(r+sz_q^i) \in q\Z$ we deduce from this equation:
$$
s[h_a + \sum_{j} z_p^jh_j + z_p^ia - z_q^ia - z_q^i z_p] \in q G.
$$
But (\ref{star}) for $p=q$ tells us that $s(h_a +\sum_j z_q^j h_j
-z_q^i z_ q) \in q G$, hence
$$
s[\sum_j(z_p^j-z_q^j)h_j + (z_p^i - z_q^i)a + z_q^i(z_q-z_p)] \in
q G.
$$
In the last expression $(z_p^i-z_q^i)a \not\in F$, thus
subtracting $(z_p^i-z_q^i)(a+z_q)\in q G$ we get
\begin{equation}\label{triplestart}
s\left[\sum_j(z_p^j-z_q^j)h_j + z_q^i(z_q-z_p) -
(z_p^i-z_q^i)z_q\right] \in q G \cap F= qF.
\end{equation}
 Now we need to look at the exact presentation of the elements $h_j$ in $F$. Let $h_j = \sum\limits_k h_j^k e_k$, where almost all coefficients $h_j^k$ are equal to $0$.
 By (\ref{psi}) and the fact that $\psi=q\psi'$ we conclude that $sh_j^i\in q \Z$ for all $j\neq i$ when restricting to the $i^{th}$-component $\Z e_i$ of $F$.
 Restricting (\ref{triplestart}) to the $i^{th}$-component of $F$ we thus obtain
$$
s[(z_p^i-z_q^i)h_i^i + (z_q^i-z_p^i) z_q^i - (z_p^i - z_q^i)
z_q^i)] \in q \Z.
$$
Therefore
$$
s(z_p^i-z_q^i)(h_i^i -  z_q^i - z_q^i) = s(z_p^i-z_q^i)(h_i^i -
2z_q^i) \in q \Z.
$$
Adding $2(z_p^i-z_q^i)(r+sz_q^i)\in q \Z$ to this expression we
obtain:
\begin{equation}\label{A}
(z_p^i- z_q^i)(sh_i^i+2r) \in q \Z.
\end{equation}
From (\ref{psiei}) we have $re_ i + sh_i + sa\in q G$, and also
$s(a+z_q)\in q G$, thus their difference
$$
(re_ i + sh_i + sa) - (sa + sz_q) = re_i + sh_i - sz_q \in q G
\cap F = q F.
$$
Restricted to the $i^{th}$-component gives:
\begin{equation}\label{B}
r+ sh_i^i + sz_q^i\in q \Z.
\end{equation}
The difference between the expression (\ref{A}) and $(z_p^i -
z_q^i)$ times the expression in (\ref{B}) gives
$$
(z_p^i - z_q^i)(r-sz_q^i) \in q \Z.
$$
Again, $(r+sz_q^i)\in q \Z$ implies
\begin{equation}\label{4r}
(z_p^i - z_q^i)(2r) \in q \Z.
\end{equation}

\noindent

Note that equation (\ref{4r}) holds for any pair of primes $p,q$
and integers $r,s$ such that $(r+sz_q^i) \in q\Z$.\\
We now distinguish two cases and make particular choices for the
data $p.q, r, s$.\\

{\bf Case $1$}: {\it There exists an index $i$ and primes $p,q$
such that $z_p^i=0$ and $z_q^i\not=0$.}\\
Note that in this case $z_q^i \not\in q\Z$ and recall that $q$
must be different from $2$. We choose integers $r,s$ relatively
prime to $q$ such that $(r+sz_q^i) \in q\Z$. Thus equation
(\ref{4r}) yields that $2rz_q^i \in q\Z$ - a
contradiction.\\

{\bf Case $2$}: {\it Not case $1$.}\\
Then for every index $i$ either $z_p^i=0$ for all primes $p$ or
$z_p^i \not=0$ for all primes $p$. In the latter we fix $i$ and claim that\\

{\it Either there exists a prime $q (\not= 2)$ such that
$(z_p^i-z_q^i)\not\in q\Z$, for some $p\in \Pi \backslash \{2\}$
or $z_p^i=z_q^i$ for all $p$ and $q$.}

\bigskip

\noindent Suppose that $(z_p^i-z_q^i)\in q\Z$ for all $p\in \Pi$.
Then obviously $(z_{p}^i-z_l^i)\in q\Z$ for all primes $p, q, l
\in \Pi$. This is impossible once we fix two primes $p$ and $l$.
Therefore, the only possibility is that $z_p^i=z_l^i$ for all
primes $p$ and $l$ in $\Pi$.\\

Now assume that there exists a prime $q (\not= 2)$ such that
$(z_p^i-z_q^i)\not\in q\Z$, for some $p\in \Pi \backslash \{2\}$.
As above we may choose integers $r,s$ relatively prime to $q$ such
that $(r+sz_q^i) \in q\Z$. Thus equation (\ref{4r}) yields that
$2r(z_p^i-z_q^i) \in q\Z$ - a contradiction.\\

The only possibility left is that for all $i$ we deduce that
$z^i_p=z^i_l$ for all $p, l\in \Pi$. However, this implies that
$z_p=z_q$ for all primes $p,q$ and hence the group $G=\langle F,
\frac{a+z}{p}:p\in \Pi\backslash\{2\}\rangle$, for a fixed element $z\in F$. In particular, one
can define a section $H\to G$ of the cellular cover $G\to H$,
given by $1\mapsto a+z$, which is not possible. This finishes the
proof of the theorem.
\end{proof}

\section{Cellular covers with cotorsion-free kernels}\label{sectionrank2}

Inspired by the previous sections we are now interested in the
following question: Can we realize every cotorsion-free abelian
group $K$ (in particular every free abelian group $K$) as the
kernel of a cellular cover of some torsion-free group of rank
two? By Theorem \ref{free-kernel} this is the best we can hope
for and by \cite[Proposition 1.4]{FGSS07} we will have to assume
that $K$ is countable since free groups are
cotorsion-free.

\subsection{Notation (see \cite{GT})} Let $R$ be a commutative ring with $1$ and a
distinguished countable multiplicatively closed subset $\bS=\{
s_n : n \in \omega \}$ such that $R$ is $\bS$-reduced and
$\bS$-torsion-free. Thus $\bS$ induces a Hausdorff topology on
$R$, taking $q_mR$ ($m \in \Z)$ as the neighborhoods of zero
where $q_m=\prod_{n < m}s_n$. We let $\Rhat$ be the $\bS$-adic
completion of $R$. We will also assume that $R$ is cotorsion-free
(with respect to $\bS$), this is to say that $\Hom(\Rhat,R) =0$.
More generally, an $R$-module $M$ is {\it $\bS$-cotorsion-free}
if $\Hom_R(\Rhat, M)=0$. We must say what it means if $M$ has
rank $\k \le\size{R}$. (Note that $R$ may not be a domain.) If
$\size{M} > \size{R}$ it suffices to let $\rk(M)=\size{M}$. If
$\size{M}\le \size{R}$, then $\rk(M)= \k$ means that there is a
free submodule $E =\bigoplus_{i < \k}Re_i$ of $M$ such that $M/E$
is $\bS$-torsion. (Note that $E$ also exists if $\size{M} >
\size{R}$.) Recall that $M$ is $\bS$-torsion if for all $m\in M$
there is $s\in \bS$ such that $sm=0$. Similarly, a submodule $N$
of $M$ is $\bS$-pure if $sM \cap N=sN$ for all $s \in \bS$. If
$M$ is $\bS$-torsion-free and $N \subseteq M$, then we denote by
$N_*$ the smallest pure submodule of $M$ containing $N$, i.e.
$N_*=\{m \in M | \exists s \in \bS \textrm{ and } sm \in N \}$.

We will write $\Hom(M,N)$ for $\Hom_R(M,N)$ and in what follows
all appearances of torsion, pure, etc. refer to $\bS$ and we will
therefore not mention the underlying set $\bS$.

\subsection{The construction}
We recently proved in \cite{GRS09} the following result:

\begin{theorem} \label{free-kernel-general}
Let $K$ be any torsion-free and reduced $R$-module of rank $\k <
\Cont$. Then there is a cotorsion-free $R$-module $G$ of rank $3$
if $\k=1$, and of rank $3\k+1$ if $2\le \k< \Cont$, with
submodule $K$ such that $\Hom(G,K)=0$ and $\Hom(G,G/K) =R \pi$
where $\pi: G\arr G/K \ (g\mapsto g+K) $ is the canonical
epimorphism. In particular,
\[ 0 \rightarrow K \rightarrow G \rightarrow G/K \rightarrow 0 \]
is a cellular cover.
\end{theorem}

If $K$ is a single copy of $R$, then the above result yields a
particular cellular cover
$$
0\to R \to G \to H\to 0
$$
where $G$ and $H$ are of rank  3 and 2, respectively. Thus $R$ is
the kernel of a cellular cover of rank two. However, if we choose
$K=R \oplus R$, then Theorem \ref{free-kernel-general} only
proves that $K$ is the kernel of a cellular cover of rank $7$. \\
We want to push this bound down to $2$ using a constructon that
goes back to A.L.S. Corner (see also \cite[Theorem 3.2]{GRS09}).
As explained before, the case of free $K$ is included and hence
we have to assume that the rank of $K$ is countable.

\begin{theorem}\label{rank2}
Let $K$ be any cotorsion-free and reduced $R$-module of finite rank $\kappa$.
Then there is a cellular exact sequence:
$$
0\to K\to G\to H
$$
where $G$ is of rank $\kappa + 2$, and therefore $H$ is of rank 2.
\end{theorem}

\begin{proof}
Let $E=\langle e_i : i \le \kappa \rangle\subseteq K$ be a
free-$R$-module, such that $K/E$ is torsion. Choose $F=\langle
f\rangle$ a free $R$-module on 1 generator. Let $C=K\oplus F$ and
$G$ be the following pure submodule of the completion $\widehat C$
of $C$:
$$
G=\langle K, F, \sum_{i=1}^\kappa w_ie_i + w f\rangle_* \subseteq
\widehat{C}
$$
where $w$ and $w_i$ ($i<\kappa$) are elements in $\widehat R$
which are algebraically independent over $K\oplus F$. For their
existence see \cite[Theorem 1.1.20]{GM90}. We claim that
\[ 0 \rightarrow K \rightarrow G \overset{\pi}{\rightarrow} H \rightarrow 0 \]
is a cellular cover, where $\pi: \widehat{K} \oplus \widehat{F}
\rightarrow \widehat{K}$ is the canonical projection. Clearly,
the rank of $G$ is then $\kappa +2$, and the cokernel $H=\pi(G)=
\langle f, wf \rangle_*$ has rank $2$.\\

We first prove that $\widehat{K} \cap G=K$ and hence $ker(\pi
\restriction_{G})=K$. Let $x \in \widehat{K} \cap G$. Then there
is $s \in \bS$ such that
\[ sx = k + f' + r(\sum\limits_{i=1}^\kappa w_ie_i + wf) \in \widehat{K} \]
for some $r \in R$ and $k \in K, f' \in F$. It follows that
\[ sx-k-r\sum\limits_{i=1}^\kappa w_ie_i = f' + rwf \] inside
$\widehat{C}$. Since $\widehat{K} \cap \widehat{F} = 0$ we
conclude that $f'+rwf = 0$ and by the algebraic independence of
$w$ also $r=0$ and $f'=0$. Thus
\[ sx-k = 0 \]
and so $sx=k \in K$ which implies $x \in K$ by the purity of $K$
in $\widehat{K}$.\\

As in the proof of Theorem 3.2 in \cite{GRS09} we need to show
that the group $G$ satisfies the desired properties, namely
$\Hom(G,K)=0$, $\Hom(G,H)= \pi R$. Therefore let $\varphi \in
\Hom(G,K)$. Thus
\[ \varphi(\sum_{i=1}^\kappa w_ie_i + w f)=\sum_{i=1}^\kappa w_i\varphi(e_i) + w
\varphi(f) \in K \] by the unique lifting of $\varphi$ to some
map from $\widehat{C}$ to $\widehat{K}$ (again denoted by
$\varphi$). Since the $w_i$ and $w$ were chosen algebraically
independent over $K$ we conclude that $\varphi(e_i)=\varphi(f)=0$
for every $i \leq \kappa$. Note that $\varphi(e_i)$ and
$\varphi(f)$ are also in $K$. Thus $\varphi \restriction_{E \oplus
F} \equiv 0$. Since the quotient $K/E$ is torsion this implies
that also $\varphi(K)=0$ by the torsion-freeness of $K$. Now, let
$x \in G$, then there is $s \in S$ such that $sx \in \left<
K,F,\sum_{i=1}^\kappa w_ie_i + w f \right>$. It follows that
$\varphi(sx)=0$ and hence also $\varphi(x)=0$ by torsion-freeness
of $K$ once more. Thus
$\varphi\equiv 0$.\\
Now, let $\varphi \in \Hom(G,H)$. As above
\[s\varphi(\sum_{i=1}^\kappa w_ie_i + w f)=\sum_{i=1}^\kappa sw_i \varphi(e_i) + sw \varphi(f) = r_ff + r_wwf  \in H \] for some $s \in S$
and coefficients $r_f,r_w \in R$. By the algebraic independence
of $w_i$ and $w$ we conclude that $\varphi(e_i)=0$ for al $i \leq
\kappa$ and hence $K$ is in the kernel of $\varphi$. Moreover,
letting $\varphi(f)=t_f f + t_wwf$ we conclude that \[ st_f w f +
st_w w^2 f = r_f f + r_w w f \] and thus $t_w=0$ and $r_f=0$.
Hence $\varphi(f)=r f$ for some $r \in R$ and therefore $\varphi$
factors through $\pi$ and induces multiplication by $r$ on $H$.
We conclude that $\varphi= \pi r \in \pi R$ and this finishes the
proof.
\end{proof}

We conclude this paper with an open question.

\begin{problem}
Can we also get a countable free group as the kernel of a rank
two cellular cover?
\end{problem}

\noindent Lutz Str\"ungmann \\
Department of Mathematics,\\ University of Duisburg-Essen,\\
Campus Essen, 45117 Essen, Germany\\
{\small e-mail: lutz.struengmann@uni-due.de}\\
\\
Jos\'e L. Rodr\'{\i}guez\\
\'Area de Geometr\'{\i}a y Topolog\'{\i}a,\\
Facultad de Ciencias Experimentales,\\ University of Almer\'{\i}a,\\
La ca{\~n}ada de San Urbano, 04120 Almer\'{\i}a, Spain\\
{\small e-mail: jlrodri@ual.es}\\

\end{document}